\documentclass[12pt,thmsa]{article}
\usepackage{amssymb}

\begin{document}

\author{Salvador Romaguera$^{1}$, Oscar Valero$^{2}$\thanks{%
The authors thank the support of Ministry of Economy and Competitiveness of
Spain, Grant MTM2012-37894-C02-01} \\
$^{1}${\small Instituto Universitario de Matem\'{a}tica Pura y Aplicada, }\\
{\small Universitat Polit\`{e}cnica de Val\`{e}ncia, 46022 Valencia, Spain }%
\\
{\small E-mail: sromague@mat.upv.es}\\
$^{2}${\small Departamento de Ciencias Matem\'{a}ticas e Inform\'{a}tica, }\\
{\small Universidad de las Islas Baleares, 07122 Palma, Spain}\\
{\small E-mail: o.valero@uib.es}}
\title{On the structure of formal balls of the balanced quasi-metric domain of words}
\date{}
\maketitle

\begin{abstract}
In ``Denotational semantics for programming languages, balanced quasi-metrics and
fixed points'' (International Journal of Computer Mathematics 85 (2008),
623-630), J. Rodr\'{i}guez-L\'{o}pez, S. Romaguera and O. Valero introduced and studied a balanced quasi-metric on any domain of (finite and infinite) words, denoted by $q_{b}$. In this paper we show that the poset of formal balls associated to $q_{b}$ has the structure of a continuous domain.\medskip

\noindent AMS Classification: 03G10, 06A06, 54E35
\smallskip

\noindent ACM Classification: F.4.3, F.3.2

\end{abstract}

\section{Introduction and preliminaries}

Throughout this paper the symbols $\Bbb{R}^{+}$ and $\Bbb{N}$ will denote
the set of all non-negative real numbers and the set of all positive integer
numbers, respectively.

Our basic references for quasi-metric spaces are \cite{Cobzas,Kun}, for general topology it is \cite{En} and for
domain theory is \cite{Gierz}.

In 1998, Edalat and Heckmann \cite{EdHec} established an elegant connection
between the theory of metric spaces and domain theory by means of the notion
of a (closed) formal ball.

Let us recall that a formal ball for a set $X$ is simply a pair $(x,r),$
where $x\in X$ and $r\in \Bbb{R}^{+}$. The set of all formal balls for $X$
is denoted by $\mathbf{B}X.$

Edalat and Heckmann observed that, given a metric space $(X,d),$ the
relation $\sqsubseteq _{d}$ defined on $\mathbf{B}X$ as
\[
(x,r)\sqsubseteq _{d}(y,s)\Leftrightarrow d(x,y)\leq r-s,
\]
for all $(x,r),(y,s)\in \mathbf{B}X,$ is a partial order on $\mathbf{B}X.$
Thus $(\mathbf{B}X,\sqsubseteq _{d})$ is a poset.

In particular, they proved the following.\newline

\textbf{Theorem 1} (\cite{EdHec}). \textit{For a metric space }$(X,d)$%
\textit{\ the following are equivalent:}

\textit{(1) }$(X,d)$\textit{\ is complete.}

\textit{(2) }$(\mathbf{B}X,\sqsubseteq _{d})$\textit{\ is a dcpo.}

\textit{(3) }$(\mathbf{B}X,\sqsubseteq _{d})\,$\textit{is a continuous
domain.}\newline

Later on, Aliakbari et al. \cite{AHPR}, and Romaguera and Valero \cite{RoVa2}
studied the extension of Edalat-Heckmann's theory to the framework of
quasi-metric spaces.\smallskip

Let us recall that a quasi-metric space is a pair $(X,d)$ where $X$ is a set
and $d:X\times X\rightarrow \Bbb{R}^{+}$ satisfies the following conditions
for all $x,y,z\in X:$

(i) $x=y\Leftrightarrow d(x,y)=d(y,x)=0;$

(ii) $d(x,y)\leq d(x,z)+d(z,y).$

The function $d$ is said to be a quasi-metric on $X.$

If the quasi-metric $d$ satisfies for all $x,y\in X$ the condition\smallskip

(i') $x=y\Leftrightarrow d(x,y)=0,\smallskip $

\noindent then $d$ is called a $T_{1}$ quasi-metric and the pair $X,d)$ is
said to be a $T_{1}$ quasi-metric space.

If $d$ is a quasi-metric on a set $X,$ then function\textbf{\ }$d^{s}$
defined as $d^{s}(x,y)=\max \{d(x,y),d(y,x)\}$ for all $x,y\in X,$ is a
metric on $X.\smallskip $

Next we recall some notions and properties of domain theory which will
useful later on.

A partially ordered set, or poset for short, is a (non-empty) set $X$
equipped with a (partial) order $\sqsubseteq .$ It will be denoted by $%
(X,\sqsubseteq )$ or simply by $X$ if no confusion arises.

A subset $D$ of a poset $X$ is directed provided that it is non-empty and
every finite subset of $D$ has upper bound in $D.$

A poset $X$ is said to be directed complete, and is called a dcpo, if every
directed subset of $X$ has a least upper bound. The least upper bound of a
subset $D$ of $X$ is denoted by $\sqcup D$ if it exists.

Let $X$ be a poset and $x,y\in X;$ we say that $x$ is way below $y,$ in
symbols $x\ll y,$ if for each directed subset $D$ of $X$ having least upper
bound $\sqcup D$, the relation $y\sqsubseteq \sqcup D$ implies the existence
of some $u\in D$ with $x\sqsubseteq u.$

A poset $X$\ is called continuous if for each $x\in X$, the set $\Downarrow
x:=\{y\in X:y\ll x\}$\ is directed with least upper bound $x$.

A continuous poset which is also a dcpo is called a continuous domain or,
simply, a domain.\smallskip

In the sequel we shall denote by $\Sigma $ a non-empty alphabet and by $%
\Sigma ^{\infty }$ the set of all finite and infinite words (or strings) on $%
\Sigma .$ We assume that the empty word $\phi $ is an element of $\Sigma
^{\infty },$ and denote by $\sqsubseteq $ the prefix order on $\Sigma
^{\infty }.$ In particular, if $x\sqsubseteq y$ and $x\neq y,$ we write $%
x\sqsubset y.$ For each $x,y\in \Sigma ^{\infty }$ we denote by $x\sqcap y$
the longest common prefix of $x\,$and $y,$ and for each $x\in \Sigma
^{\infty }$ we denote by $\ell (x)$ the length of $x.$ In particular, $\ell
(\phi )=0.$

It is well known that $\Sigma ^{\infty }$ endowed with the prefix order has
the structure of a domain.\smallskip

Usually it is defined a distinguished complete metric $d_{B}$ on $\Sigma
^{\infty },$ the so-called Baire metric (or Baire distance), which is given
by\smallskip

$\qquad d_{B}(x,x)=0$ for all $x\in \Sigma ^{\infty },\quad $and$\quad
d_{B}(x,y)=2^{-\ell (x\sqcap y)}$ for all $x,y\in \Sigma ^{\infty }$ with $%
x\neq y.\smallskip $

\noindent (We adopt the convention that $2^{-\infty }=0).\medskip $

Observe that $(\mathbf{B\Sigma }^{\infty },\sqsubseteq _{d_{B}})$ is also a
domain by Theorem 1 above.\smallskip

Recall that the classical Baire metric (or Baire distance) provides a
suitable framework to obtain denotational models for programming languages
and parallel computation \cite{deBakker1,deBakker2,deBakker3,Kahn} as well
as to study the representation of real numbers by means of regular languages
\cite{Lecomte}.\textbf{\ }However, the Baire metric is not able to decide if
a word $x$ is a prefix of another word $y,$ or not, in general. In order to
avoid this disadvantage, some interesting and useful quasi-metric
modifications of the Baire metric has been constructed. For instance:\medskip

(A) The quasi-metric $d_{w}$ defined on $\Sigma ^{\infty }$ as (compare
\cite[etc.]{Ku2,Ma,RoSch,RoVa2})$\smallskip $

$\qquad d_{w}(x,y)=2^{-\ell (x\sqcap y)}-2^{-\ell (x)}$ for all $x,y\in
\Sigma ^{\infty }.\smallskip $

(B) The quasi-metric $d_{0}$ defined on $\Sigma ^{\infty }$as (compare
\cite[etc.]{Ku2,RRVa,RoVa2,Sm1})\smallskip

$\qquad d_{0}(x,y)=0\quad $if $x$ is a prefix of $y,\smallskip $

$\qquad d_{0}(x,y)=2^{-\ell (x\sqcap y)}\quad $otherwise.$\smallskip $

(C) The $T_{1}$ quasi-metric $q_{b}$ defined on $\Sigma ^{\infty }$as
(compare \cite{RRVa})\smallskip

$\qquad q_{b}(x,y)=2^{-\ell (x)}-2^{-\ell (y)}\quad $if $x$ is a prefix of $%
y,\smallskip $

$\qquad q_{b}(x,y)=1\quad $otherwise.$\medskip $

Observe that in Examples (A) and (B) above, the fact that a word $x$ is a
prefix of another word $y$ is equivalent to say that the distance from $x$
to $y$ is exactly zero, so this condition can be used to distinguish between
the case that $x$ is a prefix of $y$ and the remaining cases for $x,y\in
\Sigma ^{\infty }.$ Observe also that $(d_{0})^{s}$ coincides with the Baire
metric while $(d_{w})^{s}$ does not.

Nevertheless, if $x,y,z\in \Sigma ^{\infty }$ satisfy $x\sqsubset y\sqsubset
z,$ one obtains $d_{w}(x,z)=d_{w}(y,z)=d_{0}(x,z)=d_{0}(y,z)=0,$ and it is
not possible to decide which word of the two, $x$ or $y$, provides a better
approximation to $z$. The quasi-metric $q_{b}$ as constructed in (C) saves
this inconvenience because if $x\sqsubset y\sqsubset z,$ it follows that $%
\ell (x)<\ell (y)<\ell (z),$ and thus $q_{b}(y,z)<q_{b}(x,z).$ Moreover, for
$x\neq \phi ,$ $x$ is a prefix of $y$ if and only if $q_{b}(x,y)<1,$ so this
condition also allows us to distinguish between the case that $x$ is a
prefix of $y$ and the rest of cases (see \cite[Remark 3]{RRVa}). We also
point out that, contrarily to $d_{w}$ and $d_{0},$ the quasi-metric $q_{b}$
has rich topological and distance properties; in particular, it is a
balanced quasi-metric in the sense of Doitchinov \cite{Do}, and consequently
its induced topology is Hausdorff and completely regular
\cite[Theorem 1 and Remark 4]{RRVa}.

By \cite[Theorem 3.1]{RoVa} (see also \cite[p. 461]{RoVa2}), $(\mathbf{%
B\Sigma }^{\infty },\sqsubseteq _{d_{w}})$ is a domain. On the other hand,
it was shown in \cite[Example 3.1]{RoVa2} that $(\mathbf{B\Sigma }^{\infty
},\sqsubseteq _{d_{0}})$ is a domain. In the light of these results, it
seems natural to wonder if $(\mathbf{B\Sigma }^{\infty },\sqsubseteq
_{q_{b}})$ is also a domain. Here we show that, indeed, this is the case.

\section{The results}

In the rest of the paper, given a quasi-metric space $(X,d),$ the way below
relation associated to $\sqsubseteq _{d}$ will be denoted by $\ll _{d}.$%
\newline

\textbf{Lemma 1} (\cite{AHPR}).\textit{\ For any quasi-metric space }$(X,d)$
\textit{the following holds:}
\[
(x,r)\ll _{d}(y,s)\Rightarrow d(x,y)<r-s.
\]

\textbf{Lemma 2.}\textit{\ Let }$(X,d)$ \textit{be a quasi-metric space. If
there is }$(x,r)\in \mathbf{B}X$ \textit{such that }$(x,r+s)\ll _{d}(x,r)$
\textit{for all }$s>0,$ \textit{then }$\Downarrow (x,r)$ \textit{is directed
and }$(x,r)=\sqcup \Downarrow (x,r).$\newline

\textit{Proof.} Obviously $\Downarrow (x,r)\neq \varnothing .$ Now let $%
(y,s),(z,t)\in \mathbf{B}X$ such that $(y,s)\ll _{d}(x,r)$ and $(z,t)\ll
_{d}(x,r).$ By Lemma 1, $d(y,x)<s-r-\varepsilon $ and $d(z,x)<t-r-%
\varepsilon $ for some $\varepsilon >0.$ Thus $(y,s)\sqsubseteq
_{d}(x,r+\varepsilon )$ and $(z,t)\sqsubseteq _{d}(x,r+\varepsilon ).$ Since
$(x,r+\varepsilon )\in \Downarrow (x,r),$ we conclude that $\Downarrow (x,r)$
is directed.

Finally, let $(z,t)$ be an upper bound of $\Downarrow (x,r).$ In
particular, we have that $(x,r+1/n)\sqsubseteq _{d}(z,t)$ for all $n,$ so $%
d(x,z)\leq r-t+1/n$ for all $n.$ Hence $d(x,z)\leq r-t,$ i.e., $%
(x,r)\sqsubseteq _{d}(z,t).$ Consequently $(x,r)=\sqcup \Downarrow (x,r).$%
\newline

A net $(x_{\alpha })_{\alpha \in \Lambda }$ in a quasi-metric space $(X,d)$
is called left K-Cauchy \cite{Ro,Sünderhauf} (or simply, Cauchy \cite{KuSch}%
) if for each $\varepsilon >0$ there is $\alpha _{\varepsilon }\in \Lambda $
such that $d(x_{\alpha },x_{\beta })<\varepsilon $ whenever $\alpha
_{\varepsilon }\leq \alpha \leq \beta .$ The notion of a left K-Cauchy
sequence is defined in the obvious manner.

Let $(X,d)$ be a quasi-metric space. An element $x\in X$ is said to be a
Yoneda-limit of a net $(x_{\alpha })_{\alpha \in \Lambda }$ in $X$ if for
each $y\in X,$ we have $d(x,y)=\inf_{\alpha }\sup_{\beta \geq \alpha
}d(x_{\beta },y).$ Recall that the Yoneda-limit of a net is unique if it
exists.

A quasi-metric space $(X,d)$ is called Yoneda-complete if every left
K-Cauchy net in $(X,d)$ has a Yoneda-limit, and it is called sequentially
Yoneda-complete if every left K-Cauchy sequence in $(X,d)$ has a
Yoneda-limit.\newline

\textbf{Lemma 3} (\cite[Proposition 2.2]{RoVa2}).\textit{\ A }$T_{1}$
\textit{quasi-metric space is Yoneda-complete if and only if it is
sequentially Yoneda-complete.}\newline

\textbf{Proposition 1. }\textit{The quasi-metric space }$(\Sigma ^{\infty
},q_{b})$ \textit{is Yoneda-complete.}\newline

\textit{Proof.} Since $(\Sigma ^{\infty },q_{b})$ is a $T_{1}$ quasi-metric
space it suffices to show, by Lemma 3, that it is sequentially
Yoneda-complete. To this end, let $(x_{n})_{n\in \Bbb{N}}$ be a left
K-Cauchy sequence in $(\Sigma ^{\infty },q_{b}).$ Then, there is $n_{1}\in
\Bbb{N}$ such that $q_{b}(x_{n},x_{m})<1$ whenever $n_{1}\leq n\leq m.$ So, $%
x_{n}$ is a prefix of $x_{m},$ i.e., $x_{n}\sqsubseteq x_{m},$ whenever $%
n_{1}\leq n\leq m.$

Now we distinguish two cases.\smallskip

Case 1.\textit{\ }There exists $n_{0}\geq n_{1}$ such that $x_{n}=x_{n_{0}}$
for all $n\geq n_{0}.$ Then, it is clear that
\[
q_{b}(x_{n_{0}},y)=\inf_{n}\sup_{m\geq n}q_{b}(x_{m},y).
\]
for all $y\in X.\smallskip $

Case 2. For each $n\geq n_{1}$ there exists $m>n$ such that $x_{n}\sqsubset
x_{m}.$ In this case, there exists $x=\sqcup \{x_{n}:n\geq n_{1}\},$ and $%
\ell (x)=\infty .$ We shall show that $x$ is the Yoneda-limit of the
sequence $(x_{n})_{n\in \Bbb{N}}.$

Indeed, we first note that $q_{b}(x_{n},x)=2^{-\ell (x_{n})}$ for all $n\geq
n_{1},$ and hence
\[
\sup_{m\geq n}q_{b}(x_{m},x)=\sup_{m\geq n}2^{-\ell (x_{m})}=2^{-\ell
(x_{n})},
\]
whenever $n\geq n_{1}.$ Therefore
\[
\inf_{n}\sup_{m\geq n}q_{b}(x_{m},x)=\inf_{n}2^{-\ell (x_{n})}=0=q_{b}(x,x).
\]
Finally, let $y\in \Sigma ^{\infty }$ such that $y\neq x.$ Since $\ell
(x)=\infty $ it follows that $x$ is not a prefix of $y,$ and thus for each $%
n\in \Bbb{N}$ there exists $m\geq \max \{n,n_{1}\}$ such that $x_{m}$ is not
a prefix of $y,$ so $q_{b}(x_{m},y)=1.$ We conclude that
\[
\inf_{n}\sup_{m\geq n}q_{b}(x_{m},y)=1=q_{b}(x,y).
\]
This finishes the proof.\newline

\textbf{Lemma 4 }(\cite{AHPR}). \textit{Let }$(X,d)$ \textit{be a
quasi-metric space.}

(a) \textit{If }$D$ \textit{is a directed subset of }$\mathbf{B}X,$ \textit{%
then }$(y_{(y,r)})_{(y,r)\in D}$ \textit{is a left K-Cauchy net in }$(X,d).$

(b) \textit{If }$\mathbf{B}X$ \textit{is a dcpo and D is a directed subset
of }$\mathbf{B}X$ \textit{having least upper bound }$(z,s),$ \textit{then }$%
s=\inf \{r:(y,r)\in D\}$ \textit{and z is the Yoneda-limit of the net }$%
(y_{(y,r)})_{(y,r)\in D}.$

(c) \textit{If }$(X,d)$ \textit{is Yoneda-complete, the poset }$(\mathbf{B}%
X,\sqsubseteq _{d})$\textit{\ is a dcpo.}\newline

\textbf{Proposition 2. }\textit{For each }$x\in \Sigma ^{\infty }$ \textit{%
such that }$\ell (x)<\infty ,$ \textit{each }$u\in \Bbb{R}^{+}$ \textit{and
each }$v>0$\textit{,} \textit{we have}
\[
(x,u+v)\ll _{q_{b}}(x,u).
\]

\textit{Proof.} Let $x\in \Sigma ^{\infty }$ with\textit{\ }$\ell (x)<\infty
,$ $u\in \Bbb{R}^{+}$ and $v>0,$ and let $D$ be a directed subset of $(%
\mathbf{B}\Sigma ^{\infty },\sqsubseteq _{q_{b}})$ whose least upper bound $%
(z,s)$ satisfies $(x,u)\sqsubseteq _{q_{b}}(z,s).$ (The existence of least
upper bound is guaranteed by Proposition 1 and Lemma 4(c)). We shall show
that there exists $(y,r)\in D$ such that $(x,u+v)\sqsubseteq
_{q_{b}}(y,r).\smallskip $

We first note that, by Lemma 4 (a), there exists $(y_{1},r_{1})\in D$ such
that $q_{b}(y_{(y,r)},y_{(y^{\prime },r^{\prime })}^{\prime })<1$ whenever $%
(y,r),(y^{\prime },r^{\prime })\in D$ with $(y_{1},r_{1})\sqsubseteq
_{q_{b}}(y,r)\sqsubseteq _{q_{b}}(y^{\prime },r^{\prime }).$ Therefore, by
the definition of $q_{b},$ we deduce that $y_{(y,r)}$ is a prefix of $%
y_{(y^{\prime },r^{\prime })}^{\prime }$ whenever $(y_{1},r_{1})\sqsubseteq
_{q_{b}}(y,r)\sqsubseteq _{q_{b}}(y^{\prime },r^{\prime }).$

Furthermore, by Lemma 4 (b), we have $s=\inf \{r:(y,r)\in D\},$ and there
exists $(y_{0},r_{0})\in D,$ with $(y_{1},r_{1})\sqsubseteq
_{q_{b}}(y_{0},r_{0}),$ such that $y_{(y,r)}$ is a prefix of $z$ whenever $%
(y_{0},r_{0})\sqsubseteq _{q_{b}}(y,r).\smallskip $

Now we distinguish two cases.$\smallskip $

Case 1. $x$ is a prefix of $z.$ Since, by assumption, $\ell (x)<\infty ,$
there exists $(y,r)\in D$ such that $(y_{0},r_{0})\sqsubseteq _{q_{b}}(y,r),$
$r<s+v,$ and $x$ is a prefix of $y_{(y,r)}.$ Then
\[
q_{b}(x,y_{(y,r)})=2^{-\ell (x)}-2^{-\ell (y_{(y,r)})}\leq 2^{-\ell
(x)}-2^{-\ell (z)}
=q_{b}(x,z)\leq u-s<u+v-r,
\]
and hence $(x,u+v)\sqsubseteq _{q_{b}}(y,r).\smallskip $

Case 2. $x$ is not a prefix of $z.$ Since, by assumption, $(x,u)\sqsubseteq
_{q_{b}}(z,s),$ we deduce that $q_{b}(u,z)=1\leq u-s.$ Choose $(y,r)\in D$
such that $r<s+v.$ Then
\[
q_{b}(x,y_{(y,r)})\leq 1\leq u-s<u+v-r,
\]
and hence $(x,u+v)\sqsubseteq _{q_{b}}(y,r).$ The proof is complete.\newline

\textbf{Theorem.}\textit{\ The poset of formal balls }$(\mathbf{B}\Sigma
^{\infty },\sqsubseteq _{q_{b}})$ \textit{is a domain.}\newline

\textit{Proof.} From Proposition 1 and Lemma 4 (c) it follows that the poset
$(\mathbf{B}\Sigma ^{\infty },\sqsubseteq _{q_{b}})$ is a dcpo, so it is
only necessary to prove that is also a continuous poset.

To this end we distinguish two cases.\smallskip

Case 1. Let $(x,r)\in \mathbf{B}\Sigma ^{\infty }$ such that $\ell
(x)<\infty .$ By Proposition 2 and Lemma 2, $\Downarrow (x,r)$ is a directed
subset of $(\mathbf{B}\Sigma ^{\infty },\sqsubseteq _{q_{b}})$ for which $%
(x,r)$ is its least upper bound.\smallskip

Case 2. Let $(x,r)\in \mathbf{B}\Sigma ^{\infty }$ be such that $\ell
(x)=\infty .$ Choose a sequence $(x_{n})_{n\in \Bbb{N}}$ of elements of $%
\Sigma ^{\infty }$ such that $\ell (x_{n})=n,$ $x_{n}\sqsubset x_{n+1}$ and $%
x_{n}\sqsubset x$ for all $n\in \Bbb{N}.$ By Lemma 4 (a), $(x_{n})_{n\in
\Bbb{N}}$ is a left K-Cauchy sequence, of distinct elements, in $(\Sigma
^{\infty },q_{b}),$ and, by Lemma 4 (b), $x$ is its Yoneda-limit.

Similarly to the proof of Proposition 2 we shall show that $%
(x_{n},2^{-n}+r)\ll _{q_{b}}(x,r)$ for all $n\in \Bbb{N},$ which implies, in
particular, that $\Downarrow (x,r)\neq \emptyset .$

Indeed, let $D$ be a directed subset of $(\mathbf{B}\Sigma ^{\infty
},\sqsubseteq _{q_{b}})$ with least upper bound $(z,t)$ such that $%
(x,r)\sqsubseteq _{q_{b}}(z,t)$. Then $t\leq r,$ and, by Lemma 4 (b), $%
t=\inf \{s:(y,s)\in D\},$ and there exists $(y_{0},s_{0})\in D$ such that $%
y_{(y,s)}$ is a prefix of $z$ whenever $(y_{0},s_{0})\sqsubseteq
_{q_{b}}(y,s).$

If $x=z,$ from the fact that $x_{n}$ is a prefix of $x$ we deduce the
existence of some $(y,s)\in D$ such that $(y_{0},s_{0})\sqsubseteq
_{q_{b}}(y,s),$ $s<t+2^{-\ell (y_{(y,s)})},$ and $x$ is a prefix of $%
y_{(y,s)}.$ Therefore
\[
q_{b}(x_{n},y_{(y,s)})=2^{-n}-2^{-\ell (y_{(y,s)})}\leq 2^{-n}+t-s\leq
2^{-n}+r-s,
\]
so that $(x_{n},2^{-n}+r)\sqsubseteq _{q_{b}}(y,s).$

If $x\neq z$ we have $q_{b}(x,z)=1\leq r-t.$ Let $(y,s)\in D$ such that $%
s<t+2^{-n}.$ Thus
\[
q_{b}(x_{n},y_{(y,s)})\leq 1\leq r-t<r+2^{-n}-s,
\]
so that $(x_{n},2^{-n}+r)\sqsubseteq _{q_{b}}(y,s).\smallskip $

Next we show that $\Downarrow (x,r)$ is directed. Indeed, let $%
(y,s),(z,t)\in \mathbf{B}\Sigma ^{\infty }$ be such that $(y,s)\ll
_{q_{b}}(x,r)$ and $(z,t)\ll _{q_{b}}(x,r).$ Since $((x_{n},2^{-n}+r))_{n\in
\Bbb{N}}$ is an ascending sequence in $(\mathbf{B}\Sigma ^{\infty
},\sqsubseteq _{q_{b}})$ with least upper bound $(x,r),$ there exists $k\in
\Bbb{N}$ such that $(x_{k},2^{-k}+r)$ is an upper bound of $(y,s)$ and $%
(z,t).$ From the fact, proved above, that $(x_{k},2^{-k}+r)\ll _{q_{b}}(x,r),
$ we deduce that $\Downarrow (x,r)$ is directed.

Finally, let $(z,t)$ be an upper bound of $\Downarrow (x,r).$ Then $%
q_{b}(x_{n},z)\leq 2^{-n}+r-t$ for all $n\in \Bbb{N}.$ Since $%
q_{b}(x,z)=\inf_{n}\sup_{m\geq n}q_{b}(x_{m},z),$ we deduce that $%
q_{b}(x,z)\leq r-t,$ and thus $(x,r)\sqsubseteq _{q_{b}}(z,t).$ Therefore $%
(x,r)$ is the least upper bound of $\Downarrow (x,r).$

We conclude that $(\mathbf{B}\Sigma ^{\infty },\sqsubseteq _{q_{b}})$ is a
domain.\newline


\begin{thebibliography}{99}
\bibitem{AHPR}  M. Aliakbari, B. Honari, M. Pourmahdian, M.M. Rezaii, The
space of formal balls and models of quasi-metric spaces, Mathematical
Structures in Computer Science 19 (2009), 337-355.

\bibitem{deBakker1}  J.W. de Bakker, E.P. de Vink, Control Flow Semantics,
Cambridge, Massachusetts: MIT Press, 1996.

\bibitem{deBakker2}  J.W. de Bakker, E.P. de Vink, A metric approach to
control flow semantics, in: Proc. 11th Summer Conference on General Topology
and Applications, Annals of the New York Academy of Sciences 806 (1996),
11-27.

\bibitem{deBakker3}  J.W. de Bakker, E.P. de Vink, Denotational models for
programming languages: Applications of Banach's fixed point theorem.
Topology and its Applications 85 (1998), 35-52.

\bibitem{Cobzas}  S. Cobza\c{s}, Functional Analysis in Asymmetric Normed
Spaces, Birkh\"{a}user, Springer Basel, 2013.

\bibitem{Do}  D. Doitchinov, On completeness in quasi-metric spaces,
Topology and its Applications 30 (1988), 127-148.

\bibitem{EdHec}  A. Edalat, R. Heckmann, A computational model for metric
spaces, Theoretical Computer Science 193 (1998), 53-73.

\bibitem{En} R. Engelking, General Topology, Monografie Mat., Vol. 60, Polish Scientific Publishers, Warszawa, 1977.

\bibitem{Gierz}  G. Gierz, K.H. Hofmann, K. Keimel, J.D. Lawson, M. Mislove,
D.S. Scott, Continuous Lattices and Domains, Encyclopedia of Mathematics and
its Applications, Vol. 93, Cambridge University Press, 2003.

\bibitem{Kahn}  G. Kahn, The semantics of a simple language for parallel
processing, in: Proc. IFIP Congress, Stockholm, Sweden, Amsterdam: Elsevier
and North-Holland, 1974, pp. 471-475.

\bibitem{Kun}  H.P.A. K\"{u}nzi, Nonsymmetric distances and their associated
topologies: About the origins of basic ideas in the area of asymmetric
topology, in: C.E. Aull, R. Lowen (Eds.), Handbook of the History of General
Topology, Vol. 3, Kluwer, Dordrecht, 2001, pp. 853-968.

\bibitem{Ku2}  H.P.A. K\"{u}nzi, Nonsymmetric topology, Bolyai Soc. Math.
Stud. 4, Topology, Szeksz\'{a}rd, 1993. Budapest, Hungary, 1995, pp. 303-338.

\bibitem{KuSch}  H.P.A. K\"{u}nzi, M.P. Schellekens, On the Yoneda
completion of a quasi-metric spaces, Theoretical Computer Science 278
(2002), 159-194.

\bibitem{Lecomte}  P. Lecomte, M. Rigo, On the representation of real
numbers using regular languages, Theory of Computing Systems, 35 (2002),
13-38.

\bibitem{Ma}  S.G. Matthews, Partial metric topology, in: Proc. 8th Summer
Conference on General Topology and Applications, Annals of the New York
Academy of Sciences 728 (1994), 183-197.

\bibitem{RRVa}  J. Rodr\'{i}guez-L\'{o}pez, S. Romaguera, O. Valero,
Denotational semantics for programming languages, balanced quasi-metrics and
fixed points, International Journal of Computer Mathematics 85 (2008),
623-630.

\bibitem{Ro}  S. Romaguera, Left K-completeness in quasi-metric spaces,
Mathematische Nachrichten 157 (1992), 15-23.

\bibitem{RoSch}  S. Romaguera, M. Schellekens, Partial metric monoids and
semivaluation spaces, Topology and its Applications 153 (2005), 948-962.

\bibitem{RoVa}  S. Romaguera, O. Valero, A quantitative computational model
for complete partial metric spaces via formal balls, Mathematical Structures
in Computer Science 19 (2009), 541-563.

\bibitem{RoVa2}  S. Romaguera, O. Valero, Domain theoretic characterisations
of quasi-metric completeness in terms of formal balls, Mathematical
Structures in Computer Science 20 (2010), 453-472.

\bibitem{Sm1}  M.B. Smyth, Quasi-uniformities: Reconciling domains with
metric spaces, in: M. Main et al. (Eds), Mathematical Foundations of
Programming Language Semantics, Third Workshop, Tulanem, 1987. Lecture Notes
in Computer Science (Berlin: Springer), Vol. 298, 1988, pp. 236-253.

\bibitem{Sünderhauf}  Ph. S\"{u}nderhauf, Quasi-uniform completeness in
terms of Cauchy nets, Acta Mathematica Hungarica 69 (1995), 47-54.
\end{thebibliography}
\end{document}